\def\titlep{Quantum mechanics and operator algebras on the Hilbert ball
(The revised)}
\documentclass[11pt]{article}
\usepackage{graphicx,ifthen}
\usepackage{amssymb}
\font\germ=eufm10 at12pt

\def\goth#1{\hbox{\germ#1}}

\newcommand{\sdag}{\scriptsize \dag}

\newcommand{\ep}{*_{\varepsilon}}
\newcommand{\vp}{\varepsilon}

\def\qed{$\square$}
\newcommand{\qedh}{\hfill\qed \\}

\newcommand{\vep}{\varepsilon}

\newtheorem{Thm}{Theorem}[section]

\newtheorem{rem}[Thm]{Remark}

\newtheorem{defi}[Thm]{Definition}
\newtheorem{lem}[Thm]{Lemma}

\newtheorem{prop}[Thm]{Proposition}
\newtheorem{prob}[Thm]{Problem}
\newtheorem{cor}[Thm]{Corollary}

\def\cal#1{\mathcal #1}

\def\pr{{\bf Proof.}\quad}

\addtocounter{footnote}{1}
\def\sftt#1{
\setcounter{equation}{0}
\addtocounter{footnote}{1}
\section{#1}
}
\def\ssft#1{\subsection{#1}}

\begin{document}
%
%
\def\autherp{Katsunori Kawamura}
\def\emailp{e-mail: kawamura@kurims.kyoto-u.ac.jp.}
\def\addressp{College of Science and Engineering Ritsumeikan University,\\
1-1-1 Noji Higashi, Kusatsu, Shiga 525-8577, Japan
}
\newcommand{\mline}{\noindent
\thicklines
\setlength{\unitlength}{.1mm}
\begin{picture}(1000,5)
\put(0,0){\line(1,0){1250}}
\end{picture}
\par
 }
\def\sd#1{#1^{\sdag}}
\def\bx{\mbox{\boldmath$x$}}
\def\bh{B_{{\cal H}}}
%
%
%
\setcounter{section}{0}
\setcounter{footnote}{0}
\setcounter{page}{1}
\pagestyle{plain}
%
%
\title{\titlep}
\author{\autherp\thanks{\emailp}
\\ 
\\
\addressp
}
\date{}
\maketitle
%
%
\begin{abstract}
Cirelli, Mani\`{a} and Pizzocchero generalized
quantum mechanics by K\"{a}hler geometry.
Furthermore they proved that any unital C$^{*}$-algebra is represented
as a function algebra on the set of pure states
with a noncommutative $*$-product as an application.
The ordinary quantum mechanics is regarded as a
dynamical system of the projective Hilbert space ${\cal P}({\cal H})$
of a Hilbert space ${\cal H}$.
The space ${\cal P}({\cal H})$
is an infinite dimensional K\"{a}hler manifold
of positive constant holomorphic sectional curvature.
In general, such dynamical system
is constructed for a general K\"{a}hler manifold of
nonzero constant holomorphic sectional curvature $c$.
The Hilbert ball $B_{{\cal H}}$ is defined by the open unit ball in ${\cal H}$
and it is a K\"{a}hler manifold with $c<0$.
We introduce the quantum mechanics on $B_{{\cal H}}$.
As an application, we show the structure of
the noncommutative function algebra on $B_{{\cal H}}$.
\end{abstract}


%
%
\sftt{Introduction}
\label{section:first}
We study a generalization of quantum mechanics
by using the theory of operator algebra and the K\"{a}hler geometry. 
The geometric quantization \cite{Woodhouse} is well-known as 
a geometrical approach of quantization.
In this theory,
the geometry means the symplectic structure of the phase space associated 
with the classical mechanics.
On the other hand,
our interest is a geometry of the quantum mechanics itself without considering
the classical mechanics and the phase space.
Hence we do not treat the problem of quantization in this paper.

A trial quantum system was considered 
by Cirelli, Mani\`{a} and Pizzocchero 
in order to study the quantum mechanics 
from a standpoint of geometry \cite{CMP90}.
The state space in a trial quantum system 
is a complete, connected, simply connected 
Hilbert K\"{a}hler manifold of constant holomorphic sectional curvature $c$.
If $c>0$, then the system is the ordinary quantum system
with the projective Hilbert space as the state space
(see also \cite{Heslot}).
Furthermore they showed that a unital C$^{*}$-algebra 
${\goth A}$ is represented by K\"{a}hler functions on the K\"{a}hler bundle
which is constructed by 
the set of pure states and the spectrum of ${\goth A}$
in \cite{CMP94}.
This result is an application of the case $c>0$. 

We examine {\it what happens when $c<0$}.
The study of the case  $c<0$ means
both a thought experiment of a new quantum mechanics
and a new theory of operator algebra
from a standpoint of geometry.
In stead of the projective Hilbert space,
we consider a framework of quantum mechanics on the Hilbert ball
%
%
\begin{equation}
\label{eqn:hilbertball}
B_{{\cal H}}\equiv\{z\in {\cal H}: \|z\|<1\}
\end{equation}
for a complex Hilbert space ${\cal H}$ as a trial quantum system.
The space $\bh$ is an example of K\"{a}hler manifold of negative 
constant holomorphic sectional curvature.
We describe its quantum mechanics by calculating objects 
appearing in its mechanics, that is,
\[
\mbox{
\begin{tabular}{l|l}
\hline
state space of system  &  $B_{{\cal H}}$ \\
\hline
observable &  K\"{a}hler function on $B_{{\cal H}}$ \\
\hline
symmetry & isometry of $B_{{\cal H}}$ \\
\hline
 transition probability &  distance of $B_{{\cal H}}$ \\
\hline
\end{tabular}
}
\]
We describe the K\"{a}hler algebra on $\bh$
as an operator algebra by solving equations of definition
of K\"{a}hler functions on $B_{{\cal H}}$.

In this section,
we provide rough overviews for each subject:
K\"{a}hler geometry,
geometrical quantum mechanics, Hilbert K\"{a}hler manifold $\bh$,
operator algebra, their relations and main results.

%
%
\ssft{K\"{a}hler geometry} 
\label{subsection:firstone}
We briefly review the definition of K\"{a}hler manifold \cite{KobNo}.
In this paper, any (infinite dimensional) manifold means a Hilbert manifold 
with respect to the Fr\'{e}chet derivative \cite{ACLM,CL,CLM,SS}. 
Let $M$ be a complex manifold with the almost complex structure $J$.
If a Riemannian metric $g$ on $M$ satisfies $g(Jx,Jx)=g(x,y)$,
then $g$ is called a {\it Hermitian metric} on $M$.
For a Hermitian metric $g$, 
the $2$-differential form $\omega$ defined by 
%
%
\begin{equation}
\label{eqn:kahler}
\omega(x, y)=g(Jx,y)
\end{equation}
is called the {\it K\"{a}hler form} of $M$.
If $d\omega=0$, then $M$ is called a {\it K\"{a}hler manifold}
with the {\it K\"{a}hler metric} $g$.
Let $R$ be the curvature tensor of the metric connection of $M$. 
A K\"{a}hler manifold $M$ is said to be {\it of
constant holomorphic sectional curvature} $c\in {\bf R}$
if $g_{p}(U_{p},R_{p}(V_{p},U_{p},V_{p}))=c$
for every $p\in M$ and any orthonormal basis $U_{p},V_{p}$ of 
every $J_{p}$-invariant two-plane of $T_{p}M$.
Let $C^{\infty}(M)$ be the  set of all smooth functions on $M$.
For $f\in C^{\infty}(M)$,
the {\it Hamiltonian vector field} $Idf$ \cite{Arnold} of $f$ is defined by
\[\omega(Idf,Y)=df(Y)\quad (Y\in TM).\]
The {\it Poisson bracket} $\{\cdot,\cdot\}$ of $M$ is defined by
\[\{f,l\}\equiv \omega(Idf,Idl)\quad(f,l\in C^{\infty}(M)).\]
The tangent space of $M$ is a Hilbert space with respect 
to the inner product associated with the K\"{a}hler metric.
We use the complexified tangent space
\[{\cal T}_{p}M\equiv T_{p}M\oplus \overline{T_{p}M}\]
where $\overline{T_{p}M}$ is the conjugate Hilbert space of $T_{p}M$.
For $f\in C^{\infty}(M)$,
the holomorphic vector fields ${\rm grad}f$ and ${\rm sgrad}f$ are
defined by
%
%
\begin{equation}
\label{eqn:gradient}
 g({\rm grad}f, \bar{Y})=\bar{\partial} f(\bar{Y}),\hspace{.2in}
\omega({\rm sgrad}f, \bar{Y})= \bar{\partial} f(\bar{Y})
\end{equation}
for any antiholomorphic vector field $\bar{Y}$ on $M$ where
$\bar{\partial}f$ is the antiholomorphic differential of $f$. 
We call ${\rm grad}f$ and ${\rm sgrad}f$ by
the {\it holomorphic gradient} and 
the {\it holomorphic skew-gradient} of $f$, respectively.
The differential $df$ of a function $f$ is written by
the holomorphic part and antiholomorphic part $\partial f +\bar{\partial}f$. 
Then $Idf={\rm sgrad}f+\overline{{\rm sgrad}f}$
and ${\rm grad}f=J{\rm sgrad}f$.

By the Kuiper's theorem \cite{Kuiper},
the tangent bundle of any infinite dimensional Hilbert manifold 
is a trivial bundle.
Furthermore if two infinite dimensional Hilbert manifolds
are homotopy equivalent, then they are diffeomorphic.
However, the K\"{a}hler structure of them are different in general.
In fact,
we treat two infinite dimensional
Hilbert K\"{a}hler manifolds of positive and negative
holomorphic sectional curvatures, respectively.

%
%
\ssft{Reformulation of quantum mechanics by K\"{a}hler geometry} 
\label{subsection:firsttwo}
We review the foundation of quantum mechanics \cite{LL,Sakurai}.
In quantum mechanics, a state is represented by 
a non zero vector in a complex Hilbert space ${\cal H}$.
Two vectors $v$ and $w$ are identified as a state 
when there exits $c\in {\bf C}$ such that $w=cv$.
The equivalence class of non zero 
vectors in ${\cal H}$ by this identification is called a {\it unit ray}
\cite{Woodhouse}.
Hence states of the system
are represented as the set of unit rays 
of ${\cal H}$.
Therefore the set of all states is the {\it projective Hilbert space}
%
%
\begin{equation}
\label{eqn:projective}
{\cal P}({\cal H})\equiv ({\cal H}\setminus\{0\})/{\bf C}^{\times}
\end{equation}
of ${\cal H}$
where ${\bf C}^{\times}\equiv \{z\in {\bf C}:z\ne 0\}$.
It is known that
${\cal P}({\cal H})$ is a K\"{a}hler manifold as
a Hilbert manifold for ${\cal H}$ with any dimension. 

In \cite{CMP90}, Cirelli, Mani\`{a} and Pizzocchero 
introduced a geometric framework of quantum mechanics.
We briefly review it here.
The operator theory of quantum mechanics is reformulated by 
the K\"{a}hler geometry of ${\cal P}({\cal H})$. 
The Schr\"{o}dinger equation is written by the equation of flow 
of a Hamiltonian vector field on ${\cal P}({\cal H})$. 
The transition probability of two states is 
the distance of two points in ${\cal P}({\cal H})$ 
with respect to the K\"{a}hler metric.
Observables are functions on ${\cal P}({\cal H})$
and the product among operators is rewritten by the $*$-product among functions
associated with the K\"{a}hler form.
A symmetry of the system is given by an isometry of ${\cal P}({\cal H})$
with respect to the K\"{a}hler metric.
\\

\noindent
\begin{tabular}{l|l|l}
\hline
{\bf physics} & {\bf operator theory} & {\bf K\"{a}hler geometry} \\
\hline
\begin{minipage}[c]{1in}
state space 
of system
\end{minipage}
 &     set of unit rays of ${\cal H}$       
 &      ${\cal P}({\cal H}) $    \\
\hline  
\begin{minipage}[c]{1.3in}
observables 
\end{minipage}
&
\begin{minipage}[c]{1.3in}
linear operators
\end{minipage}
& 
\begin{minipage}[c]{1.7in}
K\"{a}hler functions on ${\cal P}({\cal H})$
\end{minipage}
\\
\hline
symmetry &  unitary operator & isometry of ${\cal P}({\cal H}) $ \\
\hline
\begin{minipage}[c]{1.3in}
 transition \\
 probability 
\end{minipage}
 & 
\begin{minipage}[c]{1.3in}
absolute value\\
of inner product
\end{minipage}
&  K\"{a}hler distance of $ {\cal P}({\cal H}) $\\
\hline
equation of motion   & Schr\"{o}dinger equation&
\begin{minipage}[c]{1.5in}
flow equation   
of Hamiltonian vector field on ${\cal P}({\cal H})$
\end{minipage}
\\
\hline
\end{tabular}
\\

\noindent
For example, if $\dim{\cal H}=2$, then
we can identify
${\cal P}({\cal H})(= P({\bf C}^{2}))$ with the Riemann sphere $S^{2}$
as a Riemannian manifold by the stereographic projection.
Two state vectors $\psi_{1},\psi_{2}\in {\cal H}$ are corresponded to
two points $[\psi_{1}],[\psi_{2}]\in S^{2}$.
The geodesic $\ell$ through $[\psi_{1}]$ and $[\psi_{2}]$ in $S^{2}$
is always a maximal circle in $S^{2}$.

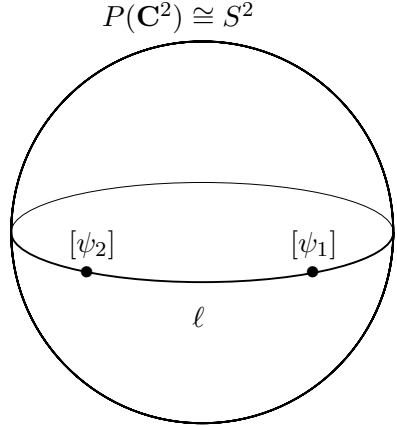
\def\sphere{
\special{pn 10}%
\special{ar 0 2260 1000 1000  0.000 19.2831853}%
\special{pn 10}%
\special{ar 0 2260 1000 261  0.0 3.14923}%
\special{pn 1}%
\special{ar 0 2260 1000 261  3.14923 7}%
}

\def\casethree{
\put(0,100){\sphere}
\put(150,-537){$\bullet$}
\put(-150,-537){$\bullet$}
\put(130,-500){$[\psi_{1}]$}
\put(-165,-500){$[\psi_{2}]$}
\put(-120,-200){$P({\bf C}^{2})\cong S^{2}$}
\put(0,-600){$\ell$}
}
\thicklines
\setlength{\unitlength}{.1mm}
\begin{picture}(1000,610)(-450,-750)
\put(40,0){\casethree}
\end{picture}

\noindent
The K\"{a}hler distance $d([\psi_{1}],[\psi_{2}])$ between
$[\psi_{1}]$ and $[\psi_{2}]$ is defined by the length 
of the geodesic between them and the following holds:
\[
d([\psi_{1}],[\psi_{2}])= \sqrt{2}\arccos |\langle \psi_{1}|\psi_{2}\rangle |
\]
where $\langle \psi_{1}|\psi_{2}\rangle$
is the inner product of $\psi_{1}$ and $\psi_{2}$ and 
we assume that $\|\psi_{1}\|=\|\psi_{2}\|=1$.

%
%
\ssft{Generalization of quantum mechanics by K\"{a}hler geometry
admitting  negative Planck constant}
\label{subsection:firstthree}
Cirelli, Mani\`{a} and Pizzocchero generalized the correspondence 
in $\S$ \ref{subsection:firsttwo}
to a general Hilbert manifold with some requirements and concluded that 
K\"{a}hler manifolds are necessary for physical ingredients
as state manifolds \cite{CMP90}.
We review the theory without the original assumption for the positivity of the
Plank constant $\hbar$.
%
%
\begin{defi}
\label{defi:gqm}
A data $(M, \omega, g, \{\Phi_{t}\}_{t\in {\bf R}})$ 
is a trial quantum system if 
the following is satisfied:
\begin{enumerate}
\item $M$ is a real smooth Hilbert manifold with a symplectic form $\omega$
 and a Riemannian metric $g$. 
\item  $\{\Phi_{t}\}_{t\in {\bf R}}$ is a 
continuous one parameter group of smooth mappings on $M$ such that
$\Phi_{t}$ preserves both $\omega$ and $g$, that is, 
$\Phi_{t}^{*}\omega =\omega$ and 
$\Phi_{t}^{*}g =g$ for each $t\in {\bf R}$.
\end{enumerate}
\end{defi}
The physical meaning of them is as follows.
The manifold $M$ is the set of pure states of a system.
The set 
%
%
\begin{equation}
\label{eqn:functions}
{\cal K}(M,{\bf R})\equiv \{
 f\in C^{\infty}(M,{\bf R});L_{Idf}g=0 \}
\end{equation}
is the set of observables where
$L_{Idf}$ is the Lie derivative by $Idf$.
The symplectic form $\omega$ 
gives a Hamiltonian system as same as classical mechanics.
The Riemannian metric $g$ gives the dispersion structure 
%
%
\begin{equation}
\label{eqn:delta}
\Delta f
\equiv \sqrt{\frac{|\hbar|}{2}g(Idf,Idf)}
\end{equation}
of an observable $f$
where the real number $\hbar$ is the Plank constant in the system.
The family $\{\Phi_{t}\}_{t\in {\bf R}}$ gives a dynamical law in the system.
%
%
%
\begin{defi}
\label{defi:axiom}
Let $(M, \omega, g, \{\Phi_{t}\}_{t\in {\bf R}})$ be 
a trial quantum system.
\begin{enumerate}
\item A symmetry of $(M, \omega, g, \{\Phi_{t}\}_{t\in {\bf R}})$
 is a smooth map which preserves both $\omega$ and $g$. 
\item 
The transition probability
between $p$ and $q$ in $M$ is the distance between
$p$ and $q$ with respect to $g$.
\end{enumerate}
\end{defi}
If the distance of $M$ is not bounded,
then we can not regard the distance as a probability.

We call $(M,\omega, g)$  the {\it state manifold} of 
$(M, \omega, g, \{\Phi_{t}\}_{t\in {\bf R}})$.
Let $T_{p}M$ and $T_{p}^{*}M$ be the tangent space 
and the cotangent space of $M$ at $p\in M$, respectively.
%
%
\begin{defi}
\label{defi:uncertainty}
\begin{enumerate}
\item The set ${\cal K}(M,{\bf R})$ in (\ref{eqn:functions}) is full
if the linear span of $\{ d_{p}f : f\in {\cal K}(M,{\bf R})\}$
 equals to $T_{p}^{*}M$ for any $p\in M$.
\item 
The manifold $M$ satisfies the uncertainty principle
if $M$ satisfies the following conditions
for $\Delta f$ in (\ref{eqn:delta}): 
\begin{enumerate}
\item 
$\Delta_{p}f\cdot \Delta_{p} l\geq 
|\frac{\hbar}{2}\{f,l\}_{p}|$ for any
 $f,l\in {\cal K}(M,{\bf R})$ and $p\in M$ where
 $\{\cdot,\cdot\}$ is the Poisson bracket with respect to $\omega$.
\item 
$\Delta_{p}f={\rm inf}\{ \lambda\geq 0 : \lambda\cdot \Delta_{p}l
\geq |\frac{\hbar}{2}\{f,l\}_{p}|, $ for every  $l\in {\cal K}(M,{\bf R})\}$
for every $f\in {\cal K}(M,{\bf R})$ and $p\in M$.
\end{enumerate}
\end{enumerate}
\end{defi}
For $p\in M$, define the linear map 
$J_{p}$ from $T_{p}M$ to $T_{p}M$
by $\omega_{p}(u,v)=g_{p}(J_{p}u,v)$ for $u,v\in T_{p}(M)$.
For $\nu\in {\bf R}\setminus\{0\}$, define
%
%
\begin{equation}\label{eqn:spd}
f*_{\nu}l\equiv f\cdot l + \frac{1}{2}\nu \left(g(Idf,Idl)
+\sqrt{-1}\omega(Idf,Idl)\right)\quad(f,l\in {\cal K}(M,{\bf C}))
\end{equation}
where ${\cal K}(M,{\bf C})\equiv
\{f:M\to {\bf C}; {\rm  Re}f,\,{\rm  Im}f\in {\cal K}(M,{\bf R})\}$.
The set ${\cal K}(M,{\bf C})$ is not closed with respect to $*_{\nu}$ in general.
But we call it the {\it $*$-product} of ${\cal K}(M,{\bf C})$ from here. 
By (\ref{eqn:spd}),
%
%
\begin{equation}
\label{eqn:hbar}
f*_{\nu}l-l*_{\nu}f=\sqrt{-1}\nu \{f,l\}.
\end{equation}
%
%
\begin{prop}\label{prop:kal}
Let $M$ be a state manifold.
\begin{enumerate}
\item 
If ${\cal K}(M,{\bf R})$ is full, then
$M$ satisfies the uncertainty principle
if and only if $J^{2}=-1$.
\item
If $J$ is integrable, that is, $\nabla J=0$, then
 for any $\nu\in {\bf R}\setminus\{0\}$,
%
%
\begin{equation}
\label{eqn:star}
f*_{\nu}(l*_{\nu}m)-(f*_{\nu}l)*_{\nu}m=0
\quad(f, l,m\in {\cal K}(M,{\bf C})).
\end{equation}
\item If ${\cal K}(M,{\bf R})$ is full and (\ref{eqn:star}) holds 
for any $ f, l,m\in {\cal K}(M,{\bf C})$, 
then $J$ is integrable.
\end{enumerate}
\end{prop}
\pr
For (i), see Proposition 4.5 in \cite{CMP90}, p 2896.
For (ii) and (iii),  see  Proposition 4.6 in \cite{CMP90}, p 2897.
\qedh 

From here, we assume that $(M,J,g)$ is a K\"{a}hler manifold. 
Then a symmetry of $(M,J,g)$ is
a K\"{a}hler (generally not holomorphic) isometry on $(M,J,g)$
and (\ref{eqn:spd}) is  rewritten as follows:
%
%
\begin{equation}
\label{eqn:hgrad}
f*_{\nu}l=f\cdot l + \nu \cdot \partial f({\rm grad}l).
\end{equation}
%
%
\begin{prop}
\label{prop:holomorphic}
\begin{enumerate}
\item If $M$ is a K\"{a}hler manifold of constant holomorphic sectional
curvature $2/\nu$, then ${\cal K}(M,{\bf C})$ is closed with respect to
the $*$-product $*_{\nu}$ in (\ref{eqn:hgrad})
and $*_{\nu}$ is associative.
\item If ${\cal K}(M,{\bf R})$ is 
full and ${\cal K}(M,{\bf C})$ is closed with respect 
to the the $*$-product $*_{\nu}$, then 
$M$ has the constant holomorphic sectional curvature $2/\nu$.
\end{enumerate}
\end{prop}
%
%
\pr
These are shown in Proposition 4.3 in the part II of \cite{CMP90}.
\qedh
From Proposition \ref{prop:holomorphic},
if $(M,J,g)$ is a K\"{a}hler manifold of non zero
constant holomorphic sectional curvature $c$,
then ${\cal K}(M,{\bf C})$ is closed with respect 
to the product $*$ defined by
%
%
\begin{equation}
\label{eqn:startwo}
f*l\equiv f\cdot l + \frac{2}{c} \cdot \partial f({\rm grad}l)\quad
(f,l\in {\cal K}(M,{\bf C})).
\end{equation}

\noindent
From (\ref{eqn:hbar}) and (\ref{eqn:startwo}),
we obtain the following relation between 
the holomorphic sectional curvature and the Planck constant: 
%
%
\begin{equation}
\label{eqn:plank}
c=\frac{2}{\hbar}.
\end{equation}

\noindent
Physical interpretations of these propositions are given in \cite{CMP90}.

%
%
\subsection{K\"{a}hler bundle and K\"{a}hler algebra}
\label{subsection:firstfour}
We review the functional representation of C$^{*}$-algebra by \cite{CMP94}.
We start from their geometric characterization of the set of pure states 
of a C$^{*}$-algebra.
%
%
\begin{defi}
\label{defi:state}
A triplet $({\cal P},p, B)$ is a K\"{a}hler bundle
if ${\cal P}$ and $B$ are topological spaces,
$p$ is a continuous surjection from ${\cal P}$ to $B$,
and 
for each $b\in B$, its fiber ${\cal P}_{b}\equiv p^{-1}(b)$
is a K\"{a}hler manifold with respect to the relative topology.
\end{defi}
We do not assume the local triviality of a K\"{a}hler bundle.
In this paper, we suppose that any manifold is a 
(possibly uncountable infinite dimensional) Hilbert manifold \cite{Schw}.
Examples of infinite dimensional K\"{a}hler Hilbert manifold is 
a projective Hilbert space, a Hilbert ball and a Loop groups \cite{Loop}.

A triplet $({\cal P},p, B)$ is a {\it uniform K\"{a}hler bundle}
 if the topology of ${\cal P}$ is a uniform topology \cite{Bourbaki}.
A triplet $({\cal P}, p,B)$ is a {\it projective K\"{a}hler bundle}
if each fiber is a projective Hilbert space.
A triplet $({\cal P}, p,B)$ is a {\it hyperbolic K\"{a}hler bundle}
if each fiber is a Hilbert ball.
A triplet $({\cal P},p, B)$ is a {\it regular state form} 
if each fiber is  a projective Hilbert space or a Hilbert ball.
%
%
\begin{defi}\label{defi:d22}
Let  $({\cal P}, p,B)$ be a K\"{a}hler bundle.
A function $f\in C^{\infty}({\cal P})$ is a K\"{a}hler function 
if $D^{2}f=0$ and $\bar{D}^{2}f=0$ where $D$ and $\bar{D}$
are the holomorphic and antiholomorphic part of 
fiberwise covariant derivative respectively.
We write ${\cal K}({\cal P})$ 
the set of all K\"{a}hler functions on $({\cal P}, p,B)$.
\end{defi}

\noindent
From Proposition \ref{prop:holomorphic}, the following holds.
%
\begin{Thm}\label{Thm:k1}
Let $({\cal P},p,B)$ be a regular state form. 
Then the set ${\cal K}({\cal P})$ of all K\"{a}hler functions is
a $*$-algebra with respect to 
the complex conjugation $f^{*}\equiv \bar{f}$ and
the $*$-product defined by
%
%
\begin{equation}
\label{eqn:productone}
(f*l)(x)\equiv f(x)\cdot l(x)+ 2\lambda_{x}\cdot\partial_{x} f({\rm grad}_{x}l)
\quad(x\in{\cal P},\,
f,l\in{\cal K}({\cal P}))
\end{equation}
where  
$\partial_{x} f$ is the holomorphic differential of $f$,
${\rm grad}_{x}l$ is the holomorphic part of gradient of $l$ with respect to
the fiberwise K\"{a}hler metric and $\lambda_{x}$ is the inverse 
of the holomorphic sectional curvature 
of the K\"{a}hler manifold of the fiber at $x$.
\end{Thm}

\noindent
In Theorem \ref{Thm:k1},
$\lambda_{x}$ is a non zero constant on $p^{-1}(b)$ for $b\in B$.
By using these preparation,
we show a geometric characterization of the set of pure states and
a functional representation of noncommutative C$^{*}$-algebras.
%
%
\begin{Thm}
\label{Thm:CMP94} 
\begin{enumerate}
\item For a unital C$^{*}$-algebra ${\goth A}$,  
the set ${\cal P}\equiv {\rm Pure}{\goth A}$ of all 
pure states of ${\goth A}$ and the spectrum $B$ of ${\goth A}$
(that is, the set ${\rm Spec}{\goth A}$ of all unitary equivalence classes 
of irreducible representations
 of ${\goth A}$),
the  natural projection $p$ from ${\cal P}$ onto $B$ is
a uniform projective K\"{a}hler bundle (=UPKB).
\item The set of all uniform continuous K\"{a}hler functions on 
a UPKB $({\cal P},p, B)$ in (i) is a unital C$^{*}$-algebra
and it is $*$-isomorphic to ${\goth A}$.
\item 
For two unital C$^{*}$-algebras ${\goth A}$ and ${\goth B}$,
their UPKB's are isomorphic if and only if 
${\goth A}$ and ${\goth B}$ are $*$-isomorphic.
\item If ${\goth A}$ is commutative, then the correspondence
in (ii) is the Gel'fand representation of commutative C$^{*}$-algebras.
\end{enumerate}
\end{Thm}

\noindent
We illustrate the uniform K\"{a}hler bundle of a C$^{*}$-algebra ${\goth A}$
as follows:

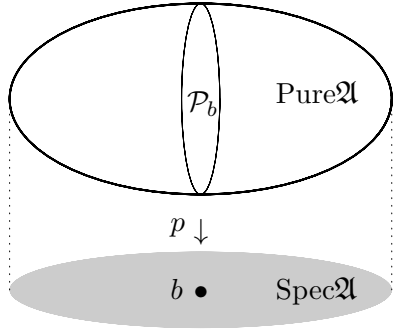
\def\sphere{
\special{pn 10}%
\special{ar -50 2260 1000 500  0.000 19.2831853}%
\special{pn 1}%
\special{ar -50 2260 100 500  0 19}%
\special{sh .2}%
\special{ia -50 3260 1000 200  0.000 19.2831853}%
\special{pn 5}%
\special{pa -1050 2260}%
\special{pa -1050 3260}%
\special{dt .04}%
\special{pn 5}%
\special{pa 950 2260}%
\special{pa 950 3260}%
\special{dt .04}%
}

\def\casethree{
\put(0,100){\sphere}
\put(100,-740){${\rm Spec}{\goth A}$}
\put(100,-480){${\rm Pure}{\goth A}$}
\put(-9,-740){$\bullet$}
\put(-40,-740){$b$}
\put(-18,-490){${\cal P}_{b}$}
\put(-40,-650){$p$}
\put(-9,-660){$\downarrow$}
}

\thicklines
\setlength{\unitlength}{.1mm}
\begin{picture}(1000,510)(-450,-820)
\put(40,0){\casethree}
\end{picture}

\noindent
By Theorem \ref{Thm:CMP94}, the category of C$^{*}$-algebras
is embedded to the category of K\"{a}hler bundles.

For each uniform continuous K\"{a}hler 
function $f$ over the K\"{a}hler bundle related a C$^{*}$-algebra ${\goth A}$,
there exists $A\in {\goth A}$ such that 
\[f(\rho)=\rho(A)\quad(\rho\in {\cal P}).\]
The norm of a function $f$ is given by 
$ \|f\|\equiv\sup_{\rho\in {\cal P}}
|(\bar{f}*f)(\rho)|^{1/2}$. 
We can regard a C$^{*}$-algebra
as a special K\"{a}hler algebra, that is,
the uniform continuous K\"{a}hler algebra
over a uniform projective K\"{a}hler bundle. 
Furthermore, the K\"{a}hler bundle
can be regard as the geometric aspects
of the noncommutative algebra such the C$^{*}$-algebra.
We showed applications of Theorem \ref{Thm:CMP94} in \cite{TD,SS}.
For the other characterization of the state space
of a C$^{*}$-algebra, see \cite{AS}.

On the other hand,
Theorem \ref{Thm:k1} holds not only projective K\"{a}hler bundles
but also hyperbolic K\"{a}hler bundles.
Therefore we have a natural question as follows.
%
%
\begin{prob}
{\it What is the K\"{a}hler algebra on a hyperbolic K\"{a}hler bundle?
How is it similar to a C$^{*}$-algebra?
}
\end{prob}
%
%
%
\ssft{K\"{a}hler geometry of the Hilbert ball}
\label{subsection:firstfive}	
We introduce the K\"{a}hler geometry of the Hilbert ball $B_{{\cal H}}$ 
in (\ref{eqn:hilbertball}) in order to consider the quantum mechanics
and operator algebra on $B_{{\cal H}}$.
Almost statements are shown without difficulty
by the analogy of finite dimensional case.
Hence their proofs are shown in Appendix \ref{section:appone}.

Let ${\cal H}$ be a complex Hilbert space 
with the inner product $\langle \cdot|\cdot\rangle $
and let $\overline{{\cal H}}$ be the conjugate Hilbert space of ${\cal H}$.
For $v\in {\cal H}$,
we define $\bar{v}\in \overline{{\cal H}}$
by $\bar{v}(w)\equiv \langle v|w\rangle$ for $w\in {\cal H}$.
Then $B_{{\cal H}}$ in (\ref{eqn:hilbertball}) is an open subset of ${\cal H}$.
The space ${\cal H}$ is a linear Hilbert manifold with the modeled space ${\cal H}$. 
Hence $B_{{\cal H}}$ is a sub Hilbert manifold of ${\cal H}$. 
The local coordinate of $B_{{\cal H}}$ is the only one which 
is induced by the inclusion map of $B_{{\cal H}}$ into ${\cal H}$. 
Therefore its tangent space $T_{z}B_{{\cal H}} $ at $z\in B_{{\cal H}}$
is uniquely defined by $T_{z}B_{{\cal H}}={\cal H}$. 
We treat $T_{z}B_{{\cal H}}$ and $\overline{T_{z}B_{{\cal H}}}$
as subspaces of ${\cal T}_{z}B_{{\cal H}}$ and 
identify $u\in T_{z}B_{{\cal H}}$, $\bar{v}\in \overline{T_{z}B_{{\cal H}}}$ 
with $(u,0)$, $(0,\bar{v})\in {\cal T}_{z}B_{{\cal H}}$, respectively.
Define the complex structure $J$ on $B_{{\cal H}}$ by
\[ J_{z} : {\cal T}_{z} B_{{\cal H}}\to {\cal T}_{z} B_{{\cal H}};
\quad J_{z}(u,\bar{v})\equiv\sqrt{-1}(u,-\bar{v})\]
for $(u,\bar{v})\in {\cal T}_{z} B_{{\cal H}}$ at $z\in B_{{\cal H}}$.

Define the Bergman type metric $g$ of $B_{{\cal H}}$ 
on ${\cal T}_{z}B_{{\cal H}} $ by
%
%
\begin{equation}\label{eqn:metric}
g_{z}((u, \bar{v}), (u^{'}, \bar{v}{'})  ) \equiv 
 k_{z}(\langle v|u^{'}\rangle + \langle v^{'}|u\rangle)
+   k_{z}^{2}(\langle v|z\rangle \langle z|u^{'}\rangle 
+\langle v^{'}|z\rangle \langle z|u\rangle)
\end{equation}
for $(u, \bar{v}), (u^{'}, \bar{v}{'}) \in {\cal T}_{z} B_{{\cal H}}$
where $k_{z}\equiv 1/(1-\|z\|^{2})$.
Then $g$ is invariant under $J$ and symmetric.
We write $g_{z}(u, \bar{v} )=g_{z}((u, 0), (0, \bar{v})  ) $. 
By this metric, $(B_{{\cal H}}, g)$ 
is a K\"{a}hler manifold of constant holomorphic sectional curvature $-2$. 
For this metric $g$,
define the K\"{a}hler form $\omega$ on $B_{{\cal H}}$ by (\ref{eqn:kahler}).

We show isometries of the K\"{a}hler manifold $\bh$.
For this purpose, we define several notations. 
Let ${\cal H}$ be a complex Hilbert space 
with the inner product $\langle\cdot|\cdot\rangle$.
We write the norm $\|\psi\|\equiv\sqrt{\langle\psi|\psi  \rangle}$
of $\psi\in {\cal H}$.
Let ${\cal L}({\cal H})$ be the set 
of all bounded linear operators on ${\cal H}$
and let $\bar{{\cal H}}$ be the conjugate Hilbert space of ${\cal H}$.
Define the new Hilbert space $\widetilde{{\cal H}}\equiv {\cal H}\oplus{\bf C}$.
We identify ${\cal H}$ as a subspace of $\tilde{{\cal H}}$.
Then any element of ${\cal L}(\tilde{{\cal H}})$ is written as follows:
%
%
\begin{equation}
\label{eqn:mat}
\left(
\begin{array}{cc}
A &x \\
\bar{y}& a
\end{array}
\right)
\end{equation}
for $A\in {\cal L}({\cal H})$,
$a\in {\bf C}$ and $x,y\in {\cal H}$
where 
$a\in {\cal L}({\bf C},{\bf C})$ is defined by $a\cdot c\equiv ac$
for $c\in {\bf C}$,
$x\in {\cal L}({\bf C},{\cal H})$ is defined by $x(c)=c\cdot x$ and
$\bar{y}\in \bar{{\cal H}}$ is defined by
$\bar{y}(\xi)\equiv \langle y|\xi\rangle$ for each $\xi\in {\cal H}$.

Define the {\it inhomogeneous unitary group} $U_{1}({\cal H})$ by
the subset of bounded linear operators $T$ on $\tilde{{\cal H}}$ satisfying 
\[T^{*}\vp T=\vp\]
where $\vp$ is the matrix 
defined on $\widetilde{{\cal H}}$ by
%
%
\begin{equation}
\label{eqn:hilbert}
\vp\equiv\left(
\begin{array}{cc}
-I_{1} & 0\\
0 & 1
\end{array}
\right)
\end{equation}
and $I_{1}$ is the identity operator on ${\cal H}$.
Then any element in $U_{1}({\cal H})$
with the form (\ref{eqn:mat})
satisfies the following conditions:
%
%
\begin{equation}
\label{eqn:conditions}
A^{*}A-\|y\|^{2}\cdot E_{y}=I_{1}, \quad
\|x\|^{2}-|a|^{2}=-1, \quad A^{*}x=ay
\end{equation}
where $E_{y}$ is the projection from ${\cal H}$ to the subspace ${\bf C}y$.
%
The action $U_{1}({\cal H})$ on $\bh$ is given by a 
{\it generalized linear fractional transformation} 
(or a {\it generalized M\"{o}bius transformation})
%
%
\begin{equation}
\label{eqn:linear}
\phi_{T}(z)\equiv \frac{Az+x}{\langle y|z\rangle+a}
\quad(z\in\bh)
\end{equation}
for $T\in U_{1}({\cal H})$ with the form (\ref{eqn:mat}).
We see that the group $U_{1}({\cal H})$ acts on $B_{{\cal H}}$ transitively.
	
%
%
\begin{Thm}\label{Thm:tr}
For $u,v\in B_{{\cal H}}$, the distance $d(u,v)$ between $u$ and $v$ 
is given by
\[ d(u,v)=\frac{1}{2}\log\frac{1-p+\sqrt{\|u-v\|^{2}-q^{2}}}
{1-p-\sqrt{\|u-v\|^{2}-q^{2}}} \]
where $p\equiv {\rm Re}\langle u|v\rangle$ 
and $q\equiv\sqrt{\|u\|^{2}\|v\|^{2}-p^{2}}$.
\end{Thm}

\noindent
The distance $d$ in Theorem \ref{Thm:tr} 
is called the {\it Poincar\'{e} distance} \cite{Kobayashi}.
%
%
\begin{cor}\label{cor:dis1}
For $u,v\in B_{{\cal H}}$, the distance $d(u,v)$ between
$u$ and $v$ is given by
\[ \tanh d(u,v)=\frac{\sqrt{\|u-v\|^{2}
-\|u\|^{2}\|v\|^{2}+({\rm Re}
\langle u|v\rangle)^{2}}}{1-{\rm Re}
\langle u|v\rangle  }\]
where the hyperbolic tangent ${\rm tanh}x$
is defined by
\[ \tanh x\equiv (e^{x}-e^{-x})(e^{x}+e^{-x})^{-1}\quad(x\in {\bf R}). \]
Especially,  ${\rm tanh}d(u,0)=\|u\|$.
\end{cor}

%
%
\ssft{Main theorem}
\label{subsection:firstsix}
We show our main statements here.
Let $\bh$ be as in (\ref{eqn:hilbertball}).
%
%
\begin{Thm}
\label{Thm:symmetryone}
\begin{enumerate}
\item
A symmetry of the system on $\bh$ is given by a generalized
linear fractional transformation in (\ref{eqn:linear})
or the mirror transformation
\[F_{W}z\equiv E_{W}z-E_{W^{\perp}}z\]
where $W$ is a closed real linear subspace of ${\cal H}$,
$E_{W}$ is the range projection of $W$ from ${\cal H}$ 
and $W^{\perp}$ is the orthogonal complementary subspace of $W$ with
respect to the real inner product ${\rm Re}\langle \cdot|\cdot\rangle$
of ${\cal H}$.
\item
For two states $\phi, \psi\in B_{{\cal H}}$, the transition probability 
$T(\phi,\psi)$ between $\phi$ and $\psi$ is given by
\[T(\phi,\psi)=\frac{1}{2}\log\frac{1-p+\sqrt{\|\phi-\psi\|^{2}-q^{2}}}
{1-p-\sqrt{\|\phi-\psi\|^{2}-q^{2}}} \]
where $p\equiv {\rm Re}\langle \phi|\psi\rangle$ and
$q\equiv\sqrt{\|\phi\|^{2}\|\psi\|^{2}-p^{2}}$.
\end{enumerate}
\end{Thm}
%
%
\pr
(i)
From Definition \ref{defi:axiom} (i) and $\S$ \ref{subsection:firstfive}, 
the statement holds.

\noindent
(ii)
From Definition \ref{defi:axiom} (ii) and Theorem \ref{Thm:tr},  
the statement holds.
\qedh

\noindent
By Theorem \ref{Thm:symmetryone} (ii),
we know that the transition probability in $B_{{\cal H}}$ depends 
on not only $|\langle \phi|\psi\rangle|$
but also the difference of phase factors and lengths of two states.

%
%
\begin{Thm}
\label{Thm:functions}
In $B_{{\cal H}}$, the set of observables is 
given as the set of smooth functions over $B_{{\cal H}}$ with the form 
%
%
\begin{equation}
\label{eqn:forms}
f_{C}(z)\equiv\frac{\langle \hat{z}|C\hat{z}\rangle }{1-\|z\|^{2}}
\end{equation}
where $C\in {\cal L}(\widetilde{{\cal H}})$
and $\hat{z}\equiv (z,1)\in  \widetilde{{\cal H}}$ for $z\in B_{{\cal H}}$.
The commutator of functions is given as follows:
\[ f_{C}*f_{C^{'}}- f_{C^{'}}*f_{C}=-\sqrt{-1}\{f_{C},f_{C^{'}}\}\]
where $\{\cdot,\cdot\}$ in
the R.H.S. is the Poisson bracket of $B_{{\cal H}}$.
\end{Thm}

We summarize our result about the trial quantum system on $\bh$ here.
The case of negative constant holomorphic curvature means
that the Planck constant in the system is {\it negative} by (\ref{eqn:plank}):
\[\hbar <0.\]
Hence this makes no sense in the usual physics.
Furthermore we do not know the affirmative experimental evidence for this system.
However, 
a physical constant often takes a unusual value in quantum field theory
in order to construct a new theory,
for example, complex mass, extra dimension and negative energy.
Therefore our thought experiment may be useful in the future
even if it is only a purely mathematical result.

By Theorem \ref{Thm:k1}, we know that 
the set ${\cal K}(B_{{\cal H}})$
of all K\"{a}hler functions on the Hilbert ball $B_{{\cal H}}$ 
is a $*$-algebra with respect to the $*$-product 
%
%
\begin{equation}
\label{eqn:hypstar}
f*l=f\cdot l -\partial f({\rm grad}l)\quad( f,l \in
{\cal K}(B_{{\cal H}}))
\end{equation}
and the complex conjugation
since $B_{{\cal H}}$ has the constant holomorphic sectional curvature $-2$.
We determine the structure of ${\cal K}(B_{{\cal H}})$.
 We define
the new product $\ep$ on the algebra ${\cal L}(\widetilde{{\cal H}})$
of all bounded linear operators on $\widetilde{{\cal H}}$ by 
\[A\ep B\equiv A\vp B \quad(A,B\in {\cal L}(\widetilde{{\cal H}})).\]
Then $({\cal L}(\widetilde{{\cal H}}), \ep)$
is a Banach $*$-algebra with the unit $\vp$ 
and the (ordinary) operator norm $\|\cdot\|$.
Then the following holds:
\[ \|A^{*}\ep I\ep A\|=\|A\|^{2}\quad (A\in {\cal L}(\widetilde{{\cal H}})).\]

Define the new norm $\|\cdot\|_{b}$ of ${\cal K}(B_{{\cal H}})$ by
\[\|f\|_{b}
=
\sup_{(z,\lambda)\in B_{{\cal H}}\times {\bf R}_{+}}
\left | t_{\lambda^{2}}(z)^{-1}\cdot ( t_{\lambda}*\bar{f}*h*f*t_{\lambda})(z)
\right|^{1/2}
\quad(f\in{\cal K}(B_{{\cal H}}))\]
where ${\bf R}_{+}\equiv \{\lambda\in {\bf R}:\lambda\geq 0\}$ and 
\[ t_{\lambda}(z)\equiv\frac{\lambda+\|z\|^{2}}{1-\|z\|^{2}},\quad
 h(z)\equiv t_{1}(z)=t_{\lambda}(z)|_{\lambda=1}\quad(
\lambda\geq 0,\,z\in B_{{\cal H}}).\]
On these preparation, 
we state the theorem of structure of the K\"{a}hler algebra 
${\cal K}(B_{{\cal H}})$ on $B_{{\cal H}}$
in Definition \ref{defi:d22}.
%
%
\begin{Thm}\label{Thm:mt}
\begin{enumerate}
\item 
The triplet $({\cal K}(B_{{\cal H}}), *, \|\cdot\|_{b})$ 
is a unital Banach $*$-algebra.
 \item ${\cal K}(B_{{\cal H}})\cong 
({\cal L}(\widetilde{{\cal H}}), \ep )$
 as a Banach $*$-algebra.
\end{enumerate}
\end{Thm}
We consider $({\cal K}(B_{{\cal H}}), *, \|\cdot\|_{b})$ 
as a deformation algebra in the theory of deformation quantization 
\cite{BFFLS, Fedosov, Kontsevich}.
Theorem \ref{Thm:mt} shows that
a deformation algebra is isomorphic to a Banach operator algebra
with a special operator product.
In general, a deformation algebra of a Poisson manifold is neither related to 
a known algebra except the Moyal algebra \cite{Moyal}, nor a Banach algebra.
Hence the statement in Theorem \ref{Thm:mt} is very rare. 
Furthermore the deformation parameter of a deformation algebra 
takes an indeterminate element in general.
On the other hand,
the deformation parameter of ${\cal K}(B_{{\cal H}})$ is 
uniquely determined by the holomorphic sectional curvature of $\bh$
by (\ref{eqn:startwo}).
In this sense, the K\"{a}hler structure of $\bh$ 
is represented by the $*$-product.
This is a new example of correspondence between geometry and algebra. 

In $\S$ \ref{section:second},
we show the proofs of statements in $\S$ \ref{subsection:firstsix}.
In $\S$ \ref{section:third},
we show examples of our results for $\bh$ when $\dim {\cal H}=1$.
In $\S$ \ref{section:fourth},
we discuss our results.

%
%
\sftt{Proofs of theorems}
\label{section:second}
We prove main theorems in $\S$ \ref{subsection:firstsix}. 
%
%
\begin{lem}\label{lem:con}
The holomorphic part of 
its Levi-Civita connection on $\bh$ is given by
\[(D_{X}Y)_{z}=k_{z}(\langle z|X\rangle Y+\langle z|Y\rangle X)
\quad(z\in B_{{\cal H}},\,X,Y\in T_{z}B_{{\cal H}}) \]
where $k_{z}=1/(1-\|z\|^{2})$.
\end{lem}
%
%
\pr
We check $D_{X}g_{z}(u,\bar{v})=0$ for $X,u\in T_{z}\bh$
and $\bar{v}\in \overline{T_{z}\bh}$ for any $z\in \bh$.
\[ 
\begin{array}{rl}
Xg_{z}(u,\bar{v}) =&\partial_{z}(g_{z}(u,\bar{v}))(X) \\
 =& k_{z}^{2}\langle v|u\rangle \langle z|X\rangle 
+ 2k_{z}^{3}\langle v|z\rangle \langle z|u\rangle \langle z|X\rangle
+ k_{z}^{2}\langle v|X\rangle \langle z|u\rangle,  \\
\\
g_{z}(D_{X}u,\bar{v})=&
k_{z}\langle v|D_{X}u\rangle +k_{z}^{2}\langle v|z\rangle \langle z|D_{X}u\rangle  \\
=& k_{z}^{2}(\langle v|u\rangle \langle z|Z|\rangle +\langle v|X\rangle \langle z|u\rangle )
+ 2k_{z}^{3}\langle v|z\rangle \langle z|u\rangle \langle z|X\rangle,  \\
\\
g_{z}(u, D_{X}\bar{v})=&0.
\end{array}
\]
Hence
$(D_{X}g_{z})(u,\bar{v})
=Xg_{z}(u,\bar{v})-g_{z}(D_{X}u,\bar{v})-g_{z}(u, D_{X}\bar{v})=0$.
By uniqueness of the Levi-Civita connection, the statement holds.
\qedh

Let ${\rm grad}f$ and ${\rm sgrad}f$ be as in (\ref{eqn:gradient}).
In \cite{Arnold,CMP94}, 
the notation $G \partial f={\rm grad}f$ and
$I\partial f={\rm sgrad}f$ are used
by defining linear mappings $G$ and $I$
on the set of vector fields on a K\"{a}hler manifold $M$.
For $M=B_{{\cal H}}$,  we show some formulae.
%
%
\begin{lem}\label{lem:grad}
Let $k_{z}\equiv 1/(1-\|z\|^{2})$.
For $z\in B_{{\cal H}}$ and $f\in C^{\infty}(B_{{\cal H}})$,
the following holds.
\begin{enumerate}
\item $ {\rm grad}_{z}f= k_{z}^{-1}\{(\bar{\partial}_{z}f)^{*}
-\langle z|(\bar{\partial}_{z}f)^{*}\rangle z\}$,
\item ${\rm sgrad}_{z}f= -\sqrt{-1}
k_{z}^{-1}\{(\bar{\partial}_{z}f)^{*}
-\langle z|(\bar{\partial}_{z}f)^{*}\rangle z\}$
\end{enumerate}
where $(\bar{\partial}_{z}f)^{*}$ is the vector satisfying 
$\langle v|(\bar{\partial}_{z}f)^{*}\rangle
=\bar{\partial}_{z}f(\bar{v})$  for any vector $v\in {\cal H}$.
\end{lem}
Lemma \ref{lem:grad} is proved in Appendix \ref{section:appone}.

We prove Theorem \ref{Thm:mt} step-by-step.
%
%
\begin{lem}\label{lem:dif1}
For $f\in C^{\infty}(B_{{\cal H}})$, 
$D^{2}f=0$ if and only if 
$\partial^{2}(k^{-1}\cdot f)=0$
where $k\in C^{\infty}(B_{{\cal H}})$ is defined by
$k(z)=(1-\|z\|^{2})^{-1}$ for $z\in B_{{\cal H}}$.
\end{lem}
%
%
\pr
By Lemma \ref{lem:con}, for $X,Y\in T_{z}B_{{\cal H}}$,
$(D^{2}_{z}f)(X,Y)=0$ if and only if
\[
\begin{array}{lll}
(\partial^{2}_{z}f)(X,Y)&=&k(z)\cdot (\partial_{z} f)(\langle z|X\rangle Y
+\langle z|Y\rangle X)\\
&=& k(z)\cdot\{\partial_{z}^{2}((1-k^{-1})\cdot f)(X,Y) 
-\|z\|^{2}\cdot \partial_{z}^{2}f(X,Y)\}.
\end{array}
\]
By multiplying $k(z)^{-1}$ at both sides,
\[
 \begin{array}{lll}
(1-\|z\|^{2})(\partial^{2}_{z}f)(X,Y)&=&\partial_{z}^{2}(
(1-k^{-1})\cdot f)(X,Y) -\|z\|^{2}\cdot \partial_{z}^{2}f(X,Y).
\end{array}
\]
It is equivalent to
$\partial^{2}_{z}(k^{-1}\cdot f)(X,Y)=0$. 
\qedh
%
%
\begin{lem}\label{lem:dif2}
If $l\in C^{\infty}(B_{{\cal H}})$ satisfies 
both $\partial^{2}l=0$ and $\bar{\partial}^{2}l=0$, then 
there exist $a\in {\bf C}$, $A\in{\cal L}({\cal H})$ and $u,v\in {\cal H}$
such that 
%
%
\begin{equation}\label{eqn:form}
 l(z)=\langle z|Az\rangle +\langle u|z\rangle+\langle z|v\rangle+a.
\end{equation}
\end{lem}
%
%
\pr
By assumption of $\partial^{2}l=0$, 
there exists $u_{0}\in {\cal H}$ such that 
%
%
\begin{equation}\label{eqn:e11}
(\partial_{z}l)(X)=\langle u_{0}|X\rangle \quad(z\in B_{{\cal H}})
\end{equation}
because $(\partial_{z}l)(X)$ is bounded linear with respect to $X$.
From this,
for $\bar{\partial}^{2}l=0$,
we see that there exists $v\in {\cal H}$ such that 
$(\bar{\partial}_{z}l)(\bar{X})=\langle X|v\rangle$ for $z\in B_{{\cal H}}$.
From (\ref{eqn:e11}), there exists $c_{0}\in{\bf C}$ such that
%
%
\begin{equation}\label{eqn:e12}
 l(z)=\langle u_{0}|z\rangle +c_{0}.
\end{equation}
Consider $u_{0}=u_{0}(\bar{z})$ and $c_{0}=c_{0}(\bar{z})$. Then
\[ (\bar{\partial}_{z}l)(\bar{X})=
\bar{\partial}_{z}\{ \langle u_{0}(z)|z\rangle +c_{0}(\bar{z}) \}(\bar{X})
=\langle \partial_{z}u_{0}(X)|z\rangle
+\bar{\partial_{z}}c_{0}(\bar{X}).\]
Because $(\partial_{z}u_{0})(X)$ and $\bar{\partial_{z}}c_{0}(\bar{X})$ 
are bounded linear with respect to $X$ and $\bar{X}$, respectively,
there exist $B\in {\cal L}({\cal H})$ and $v\in {\cal H}$ such that
$\partial_{z}u_{0}(X)=BX$ and
$\bar{\partial_{z}}c_{0}(\bar{X})=\langle X|v\rangle$.
Therefore, there exist $u\in {\cal H}$ and $a\in {\bf C}$ such that 
%
%
\begin{equation}\label{eqn:ins}
u_{0}(z)=Bz+u, \quad c_{0}(z)=\langle z|v\rangle+a. 
\end{equation}
By inserting (\ref{eqn:ins}) into (\ref{eqn:e12}) and let $A\equiv B^{*}$, 
we obtain (\ref{eqn:form}). \qedh

For $z\in B_{{\cal H}}$, define $\hat{z}\equiv (z,1)\in \tilde{{\cal H}}$.
%
%
\begin{defi} \label{defi:d1}
For $C\in {\cal L}(\tilde{{\cal H}})$ with the form (\ref{eqn:mat}),
define the function $f_{C}$ on $B_{\cal H}$ by
\[f_{C}(z)\equiv
\frac{a+\langle y|z\rangle +\langle z|x\rangle
+\langle z|Az\rangle}{1-\|z\|^{2}}.\]
\end{defi}
We see that 
$f_{C}(z)= k(z)\cdot \langle \hat{z}|C\hat{z}\rangle$
for each $z\in B_{{\cal H}}$.
We conclude the next corollary.
%
%
\begin{cor}\label{cor:c10}
If $f\in C^{\infty}(B_{{\cal H}})\cap {\cal K}(B_{{\cal H}})$,
then there exists $C\in {\cal L}(\tilde{{\cal H}})$ such that $f=f_{C}$.
\end{cor}
%
%
\begin{prop}\label{prop:p1}
For $C,C^{'}\in {\cal L}(\tilde{{\cal H}})$,
$f_{C}*f_{C^{'}}=f_{C\vp C^{'}}$.
\end{prop}
%
%
\pr
We calculate the noncommutative part 
$\partial_{z}f_{C}({\rm grad}_{z}f_{C^{'}})$
of the $*$-product $f_{C}*f_{C^{'}}$.
By using (\ref{eqn:ew2}) and (\ref{eqn:ew3}) in Lemma \ref{lem:a1},
\[
\begin{array}{rl}
\partial_{z}f_{C}({\rm grad}_{z}f_{C^{'}})=&
 \langle (\partial_{z}f_{C})^{*}|{\rm grad}_{z}f_{C^{'}}\rangle \\
=&\langle \{k_{z}E_{1}C^{*}\hat{z}
+k_{z}^{2}\langle \hat{z}|C^{*}\hat{z}\rangle z\}|\{
E_{1}C^{'}\hat{z}+\langle e_{2}|C^{'}\hat{z}\rangle z\}\rangle \\
=&k_{z}\langle E_{1}C^{*}\hat{z}|C^{'}\hat{z}\rangle 
+k_{z}^{2}\overline{\langle \hat{z}|C^{*}\hat{z}\rangle }\langle z|E_{1}C^{'}\hat{z}\rangle \\
&+k_{z}\langle E_{1}C^{*}\hat{z}|z\rangle \langle e_{2}|C^{'}\hat{z}\rangle 
+k_{z}^{2}\langle \hat{z}|C\hat{z}\rangle \langle e_{2}|C^{'}\hat{z}\rangle \|z\|^{2} \\
=& f_{C}(z)f_{C^{'}}(z)-k_{z}^{2}\langle \hat{z}|C\hat{z}\rangle \langle e_{2}|C^{'}\hat{z}\rangle \\
&+k_{z}\langle E_{1}C^{*}\hat{z}|C^{'}\hat{z}\rangle
+k_{z}\langle \hat{z}|Cz\rangle \langle e_{2}|C^{'}\hat{z}\rangle 
+k_{z}^{2}\langle \hat{z}|C\hat{z}\rangle \langle e_{2}|C^{'}\hat{z}\rangle \|z\|^{2}\\
=& f_{C}(z)f_{C^{'}}(z)-k_{z}\langle \hat{z}|Ce_{2}\rangle \langle e_{2}|C^{'}\hat{z}\rangle  
+k_{z}\langle E_{1}C^{*}\hat{z}|C^{'}\hat{z}\rangle. \\
\end{array}
\]
From this, we obtain that
%
%
\begin{equation}
\label{eqn:eqas}
\partial_{z}f_{C}({\rm grad}_{z}f_{C^{'}})=f_{C}(z)f_{C^{'}}(z)
-k_{z}\langle \hat{z}|C\{\langle e_{2}|C^{'}
\hat{z}\rangle e_{2}-E_{1}C^{'}\hat{z}\}\rangle. 
\end{equation}
Since 
\[ \langle e_{2}|C^{'}\hat{z}\rangle e_{2}=
\left(\begin{array}{cc}
0&0\\
0& 1
\end{array}
\right)
C^{'}\hat{z},\hspace{.2in}E_{1}C^{'}\hat{z}=
\left(\begin{array}{cc}
I_{1}&0\\
0& 0
\end{array}
\right)C^{'}\hat{z},\]
the R.H.S. of (\ref{eqn:eqas}) is 
$f_{C}(z)f_{C^{'}}(z)-k_{z}\langle \hat{z}|C\vp C^{'}\hat{z}\rangle  
= f_{C}(z)f_{C^{'}}(z)-f_{C\vp C^{'}}(z)$.
Therefore, 
\[(f_{C}*f_{C^{'}})(z)=f_{C}(z)f_{C^{'}}(z)
-\partial_{z}f_{C}({\rm grad}_{z}f_{C^{'}})
=f_{C\vp C^{'}}(z).\]
Hence the statement holds. \qedh

\noindent
{\bf Proof of Theorem \ref{Thm:mt}.}
In the first step, we show that
$({\cal K}(B_{{\cal H}}),*)$ and $({\cal L}(\widetilde{{\cal H}}), \ep )$
are isomorphic as a $*$-algebra.
By Definition \ref{defi:d1}, the mapping
$f$ from ${\cal L}(\widetilde{{\cal H}})$ to ${\cal K}(B_{{\cal H}})$,
$C\mapsto f_{C}$,
is defined. 
By Corollary \ref{cor:c10}, $f$ is surjective.
If $f_{C}=0$ for 
$C=\left(
\begin{array}{cc}
A & \bar{u}\\
v & a
\end{array}
\right )\in {\cal L}(\widetilde{{\cal H}})$,
then by comparing degree of $z\in B_{{\cal H}}$ in $f_{C}$,
$\langle z|Az\rangle =0, 
\,\langle u|z\rangle +\langle z|v\rangle =0,\, a=0$ for any $z\in B_{{\cal H}}$.
Hence $A=0$. Since $\langle u|z\rangle $ is affine 
and $\langle z|v\rangle $ is conjugate affine  with respect to $z$,
we obtain $u=v=0$. Therefore $C=0$.
Hence $f$ is injective.
Clearly, $f$ is linear and $f_{C^{*}}=\bar{f_{C}}$.
By Proposition \ref{prop:p1}, $f$ is a $*$-isomorphism.

By definition of the norm $\| \cdot \|_{b}$,
\[\|f_{C}\|_{b}^{2}
=\sup_{\lambda\geq 0,\,z\in B_{{\cal H}}}
\left | t_{\lambda^{2}}(z)^{-1}\cdot 
( t_{\lambda}*\bar{f_{C}}*h*f_{C}*t_{\lambda})(z)
\right|. \]
We can write $ t_{\lambda}=f_{T_{\lambda}}$ by
$ T_{\lambda}\equiv\left(
\begin{array}{cc}
I_{1} & 0\\
0 & \lambda
\end{array}
\right )$ and $h=f_{I}$.
From Proposition \ref{prop:p1} and (\ref{eqn:forms}),
\[
 \|f_{C}\|_{b}^{2}
=
\sup_{(z,\lambda)\in B_{{\cal H}}\times {\bf R}_{+} }
\left | \left(\frac{\lambda^{2}+\|z\|^{2}}{1-\|z\|^{2}} \right)^{-1} \cdot 
( f_{T_{\lambda}\ep C^{*}\ep I\ep C\ep T_{\lambda}})(z)
\right|.
\]
Since
\[
\left(\frac{\lambda^{2}+\|z\|^{2}}{1-\|z\|^{2}} \right)^{-1} \cdot 
( f_{T_{\lambda}\ep C^{*}\ep I\ep C\ep T_{\lambda}})(z)
 =
 \frac{\langle (-z,\lambda)|C^{*}C(-z,\lambda)\rangle }
{\lambda^{2}+\|z\|^{2}}=
\frac{ \|C(-z,\lambda)\|^{2} }{ \|(-z,\lambda)\|^{2} },
\]
we obtain $ \|f_{C}\|_{b}^{2}=\|C\|^{2}$.
Hence $f$ is an isometry from the Banach space 
 $({\cal L}(\widetilde{{\cal H}}), \|\cdot\|)$ to 
$({\cal K}(B_{{\cal H}}),\|\cdot\|_{b})$. 
Therefore $({\cal L}(\widetilde{{\cal H}}),*, \|\cdot\|)$
is a  Banach $*$-algebra with the unit $f_{J}\equiv 1$
because $({\cal L}(\widetilde{{\cal H}}), \ep, \|\cdot\|)$
is a Banach $*$-algebra. We finish to show (i). 
The statement (ii) is clear from the proof of (i).
\qedh

\noindent
{\bf Proof of Theorem \ref{Thm:functions}.}
An observable in the system is given by
a function $f$ on $B_{{\cal H}}$ such that $ D^{2}f=0$ and $\bar{D}^{2}f=0$
where $D,\bar{D}$ are the holomorphic, the antiholomorphic
part of covariant derivative of $B_{{\cal H}}$ respectively.
Therefore it is calculated in Theorem  \ref{Thm:mt}.
Later equation follows from Theorem \ref{Thm:k1}.
\qedh

%
%
\begin{rem}{\rm 
As same as the projective Hilbert space case,
we can consider a manifold domain \cite{CLM} of $B_{{\cal H}}$
for a unbounded selfadjoint operator $H$ with dense domain ${\cal D}
\subset {\cal H}$. For such $(H, {\cal D})$,
we can make ${\cal D}$  as a Hilbert space with 
the inner product $\langle \cdot|\cdot\rangle _{H}$
\[ \langle x|y\rangle_{H}\equiv \langle x|y\rangle+\langle Hx|Hy\rangle\quad
(x,y\in {\cal D}).\]
 Let $B_{{\cal D}}\equiv B_{{\cal H}}\cap {\cal D}$
is also a K\"{a}hler (Hilbert) manifold.
The inclusion mapping $\iota$ of $B_{{\cal D}}$ into $B_{{\cal H}}$
is an isometry and the range of its differential $d\iota$
is dense in each tangent space of $B_{{\cal H}}$.
}
\end{rem}

%
%
\begin{rem}{\rm The new algebra 
$({\cal L}(\widetilde{{\cal H}}), \ep)$ and 
the algebra $({\cal L}(\widetilde{{\cal H}}),\cdot)$ with 
the (ordinary) operator product $\cdot$ are isomorphic as
a unital Banach algebra by the mapping
\[ A\mapsto  \vp A.\] 
However, $({\cal L}(\widetilde{{\cal H}}), \ep)$ 
and $({\cal L}(\widetilde{{\cal H}}),\cdot)$ are not $*$-isomorphic
because $({\cal L}(\widetilde{{\cal H}}), \ep)$
does not satisfy the C$^{*}$-condition.
We explain this as follows: 
Let $v\in {\cal H}$ be a unit vector 
and let $E$ be the one dimensional projection from ${\cal H}$ to ${\bf C}v$.
Then for an operator 
\[A\equiv \left(
\begin{array}{cc}
E & 0\\
v & 0
\end{array}
\right),
\]
the following equations hold:
\[ A^{*}A=
\left(
\begin{array}{cc}
2E & 0\\
0 & 0
\end{array}
\right), \hspace{.3in}
A^{*}\vp A=0.
\]
Hence, $\|A\|^{2}=\|A^{*}A\|\ne\|A^{*}\vp A\|=\|A^{*}\ep A\|$.
}
\end{rem}

%
%
\begin{rem}{\rm
If we consider ${\cal K}(B_{{\cal H}})$ as a 
functional representation of some kind of operator algebra
and ${\cal K}(B_{{\cal H}})$ has a norm, 
then we will expect that ${\cal K}(B_{{\cal H}})$
becomes a norm algebra.
But if we treat ${\cal K}(B_{{\cal H}})$ as a function
algebra on $B_{{\cal H}}$, then we should
 consider a norm which is defined on only $B_{{\cal H}}$.
Therefore the norm $\|\cdot\|_{b}$ is not suitable for ${\cal K}(B_{{\cal H}})$.

For example, we define the other norm $\|\cdot\|_{d}$ by
\[\|f\|_{d}\equiv 
\sup_{z\in B_{{\cal H}}}|h(z)^{-1}\cdot (\bar{f}*h*f)(z)|^{1/2}\quad
(f\in {\cal K}(B_{{\cal H}})).\]
Then 
\[\|f_{C}\|_{d}=\|C\|_{s}\equiv 
\sup_{z\in B_{{\cal H}}}
\frac{\|C(z,1)\|}{\|(z,1)\|}\qquad
(C\in {\cal L}(\widetilde{{\cal H}})).
\]
Both $({\cal K}(B_{{\cal H}}), \|\cdot\|_{d})$ and
$({\cal L}(\widetilde{{\cal H}}), \|\cdot\|_{s})$
are Banach spaces and they are $*$-isomorphic 
  as a $*$-algebra and isometric but not Banach $*$-algebras
 because they don't satisfy the condition
\[ \|A\ep B\|_{s}\leq\|A\|_{s}\cdot\|B\|_{s}
\quad(A,B\in {\cal L}(\widetilde{{\cal H}})).\]
For example, 
let $v\in {\cal H}$, $\|v\|=1$ and 
let $E :{\cal H}\to {\bf C}v$ be the
one dimensional projection.
Define two matrices $C$ and C$^{'}$ on 
$\widetilde{{\cal H}}\equiv {\cal H}\oplus {\bf C}$ by
\[ C\equiv
\left(
\begin{array}{cc}
E &0\\
0&0
\end{array}
\right),\quad
C^{'}\equiv
\left(
\begin{array}{cc}
0 &\bar{v}\\
0&0
\end{array}
\right).
\]
Then 
$\|C\ep C^{'}\|_{s}\not \leq \|C\|_{s}\cdot\|C^{'}\|_{s}$.

}
\end{rem}

%
%
\sftt{$1$-dimensional case --- unit open disc}
\label{section:third}
We consider the quantum mechanics on $\bh$ 
for the case of $\dim_{\bf C}{\cal H}=1$.

%
%
\ssft{Geometry and K\"{a}hler algebra}
\label{subsection:thirdone}
In this case,
$\bh$ is the unit open disc $D$ in ${\bf C}$ defined by
\[D\equiv \{z\in {\bf C}:|z|<1\}.\]
The space $D$ is well-known as an example 
of hyperbolic complex space \cite{hel,Kobayashi}.
We review properties of $D$.
The K\"{a}hler metric $g$ on $D$ is given 
by the {\it Poincar\'{e} metric} defined by
\[g_{z}(\bar{v}, u)=2\bar{v}u/(1-|z|^{2})^{2}
\quad(z\in D,\, \bar{v},u\in {\bf C})\]
where $g$ is scaled such that the curvature of $(D,g)$ is $-2$.
For $z_{1},z_{2}\in D$,
the geodesic $\ell_{1}$ between $z_{1}$ and $z_{2}$ 
is the circular arcs perpendicular to the
boundary $|z|=1$ (Fig. 4.1).
In particular,
if $z_{1}$, $z_{2}$ and $0$ are on the same straight line $\ell_{2}$,
then $\ell_{2}$ is the geodesic between $z_{1}$ and $z_{2}$
(Fig. 4.2).

\def\circlefill{
\special{sh 0.2}%
\special{ia 0 2260 1000 1000  0.000 19.2831853}%
}
\def\arcone{
\special{pn 10}%
\special{ar 2500 1490 861 861  2.61000 4.220923}%
}
\def\lineone{
\special{pn 10}%
\special{pa 2500 2740}%
\special{pa 3900 1345}%
\special{fp}%
}
\def\casethree{
\put(0,100){\circlefill}
\put(-300,-190){\arcone}
\put(170,-437){$\bullet$}
\put(130,-645){$\bullet$}
\put(150,-400){$z_{1}$}
\put(95,-670){$z_{2}$}
\put(0,-300){$D$}
\put(100,-500){$\ell_{1}$}
}
\def\casefour{
\put(9,109){\circlefill}
\put(0,-437){$0$}
\put(10,-485){$\bullet$}
\put(127,-367){$\bullet$}
\put(-110,-605){$\bullet$}
\put(150,-400){$z_{1}$}
\put(-105,-650){$z_{2}$}
\put(0,-300){$D$}
\put(100,-500){$\ell_{2}$}
\put(-800,50){\lineone}
}
\thicklines
\label{fig:one}
\setlength{\unitlength}{.1mm}
\begin{picture}(1000,630)(-450,-800)
\put(-240,0){\casethree}
\put(-310,-780){{\bf Fig. 4.1}}
\put(430,-9){\casefour}
\put(390,-780){{\bf Fig. 4.2}}
\end{picture}

If $\dim_{{\bf C}} {\cal H}\ne 1$, then
the geodesic between $z_{1}$ and $z_{2}$ in $\bh$ is always 
transformed to the straight line $l$ through $0$ by a suitable isometry of $\bh$.
In this case, there exists a $1$-dimensional complex subspace $V\subset {\cal H}$ 
containing $l$ such that $\{v\in V:\|v\|<1\}$ is isometric to $D$.
Hence the above illustration also makes sense in the general case. 

The group of all isometries of $D$ is generated by
a linear fractional transformation $\phi_{X}$ by $X\in M_{2}({\bf C})$ 
with the form
%
%
\begin{equation}
\label{eqn:uone}
X=\left(
\begin{array}{cc}
s&r\\
\bar{r}& \bar{s}
\end{array}
\right),\quad |s|^{2}-|r|^{2}=1
\end{equation}
and the complex conjugation $z\mapsto \bar{z}$.

The K\"{a}hler algebra on $D$
is $*$-isomorphic to $(M_{2}({\bf C}),*_{\vep})$ where
$*_{\vep}$ is the product of $M_{2}({\bf C})$ defined by
\[A*_{\vep}B\equiv A\vep B\quad(A,B\in M_{2}({\bf C}))\quad
\vep\equiv {\rm diag}(-1,1)\in M_{2}({\bf C}).\]

%
%
\ssft{Dynamics of states in $D$}
\label{subsection:thirdtwo}
In the trial quantum system of $D$,
the dynamical law is given by
a continuous one parameter group 
$\{\Phi_{t}\}_{t\in {\bf R}}$ 
of K\"{a}hler isometries on $D$ by Definition \ref{defi:gqm}.
From (\ref{eqn:uone}), the Lie algebra ${\goth g}$ 
of the group of all isometries of $D$ is 
given as
\[\left\{\left(\begin{array}{cc}
\sqrt{-1}a&b \\
\bar{b}&-\sqrt{-1}a\\
\end{array}\right):
a\in {\bf R},\,b\in {\bf C}
\right\}.
\]
For a given $\{\Phi_{t}\}_{t\in {\bf R}}$,
there exists $X\in {\goth g}$ such that $\Phi_{t}=\phi_{\exp tX}$ for each $t$.
From this,
we can concretely describe the time evolution as
\[z(t)\equiv \phi_{\exp tX}(z)\quad (z\in D,\,t\in {\bf R}).\]
For $X\equiv 
\left(
\begin{array}{cc}
\sqrt{-1}a&b\\
\bar{b}&-\sqrt{-1}a\\
\end{array}
\right)\in {\goth g}$,
let $\alpha\equiv |b|^{2}-a^{2}$.
If $\alpha>0$, then
\[z(t)=\frac{(\sqrt{\alpha}+\sqrt{-1}a\tanh \sqrt{\alpha}t)z
+b\tanh \sqrt{\alpha}t}{(\bar{b}\tanh \sqrt{\alpha}t) z
+\sqrt{\alpha}-\sqrt{-1}a\tanh \sqrt{\alpha}t}.
\]
If $\alpha<0$, then
\[z(t)
=\frac{(\sqrt{-\alpha}+\sqrt{-1}a\tan \sqrt{-\alpha}t)z
+b\tan \sqrt{-\alpha}t}{(\bar{b}\tan \sqrt{-\alpha}t) z
+\sqrt{-\alpha}-\sqrt{-1}a\tan \sqrt{-\alpha}t}.
\]
In addition to this condition,
if $b=0$, then 
\[z(t) =e^{2\sqrt{-1}at}z.\]
If $\alpha=0$, then 
$z(t)=z$ for each $t$.
In $\S$ 2.1 of \cite{Sakurai}, 
the Schr\"{o}dinger equation is derived from
the assumption that the Hamiltonian is the generator of time evolution \cite{Goldstein}.
More properly, the Hamiltonian $\times \sqrt{-1}$
in the ordinary quantum mechanics
is chosen as the infinitesimal generator of a one-parameter unitary group
on the Hilbert space.
Therefore an element $H\in \frac{1}{\sqrt{-1}}{\goth g}$ is a candidate
for a Hamiltonian operator of the trial quantum mechanics on $D$. 

%
%
\sftt{Discussion}
\label{section:fourth}
We discuss  relations between cases of ${\cal P}({\cal H})$ and $\bh$.
We show the quantum mechanics of the Hilbert ball 
including that of the projective Hilbert space 
as a special case.

Assume that $H\in {\cal L}({\cal H})$ is selfadjoint
and $a\in {\bf R}$.
Define the operator $X$ on $\widetilde{{\cal H}}$ by
\[X\equiv \left(
\begin{array}{cc}
\sqrt{-1}H & 0 \\
0 & \sqrt{-1}a
\end{array}
\right).\]
Then $\exp tX\in U_{1}({\cal H})$ for each $t\in {\bf R}$
because it satisfies (\ref{eqn:conditions}).
Let $z(t)\equiv \phi_{\exp tX}(z)$ for $z\in \bh$.
Then $z(t)= e^{\sqrt{-1}(H-aI)t}z$.
If we replace $H$ by $-H+Ia$,
then
\[z(t)= e^{-\sqrt{-1}Ht}z.\]
This is the solution of
the Schr\"{o}dinger equation with the Hamiltonian $H$:
\[Hz(t) =\sqrt{-1}\frac{d}{dt}z(t).\] 

From Corollary \ref{cor:dis1}, 
the inner product $\langle\cdot|\cdot\rangle$ of ${\cal H}$ is recovered from
the K\"{a}hler distance $d$ of $\bh$ as follows.
If $u,v\in \bh$ satisfy $\|u\|=\|v\|$, then
\[
\begin{array}{lll}
\langle u|v\rangle &=&2\sinh^{2}(\frac{d(u,v)}{2})
+\tanh^{2}d(u,0){\rm cosh} d(u,v)\\
&&+\sqrt{-1}\{2\sinh^{2}(\frac{d(\sqrt{-1}u,v)}{2})
+\tanh^{2}d(u,0)\cosh d(\sqrt{-1}u,v)\}.
\end{array}
\]

\section*{Appendix}
\appendix

%
%
\sftt{Proof of lemmata}
\label{section:appone}
%
%
%
\begin{lem}\label{lem:16}
For $z\in B_{\cal H}$ and $0\leq t\leq 1$,
define $\theta(t)\in \bh$ by
%
%
\begin{equation}
\label{eqn:line}
\theta(t)\equiv t\cdot z.
\end{equation}
Then $\theta$ is the geodesic between $0$ and $z$. 
The length $L(\theta)$ of $\theta$ is given by
\[ L(\theta) =\frac{1}{2}\log \frac{1+\|z\|}{1-\|z\|}.\]
\end{lem}
%
%
\pr
We show $L(\theta)\leq L(\psi)$ for any curve $\psi$ between $0$ and $z$.
By definition of length of curve
%
%
\begin{equation}
\label{eqn:ie}
L(\psi)=\int_{0}^{1}\sqrt{K(t)}\,dt
\end{equation}
where $K(t)\equiv g_{\psi(t)}(\dot{\psi}(t),\dot{\psi}(t))$.
Then
\[K(t)=\frac{\|\dot{\psi}(t)\|^{2}}{1-\|\psi(t)\|^{2}}
-\frac{|\langle \dot{\psi}(t)|\psi(t)\rangle|^{2}}{(1-\|\psi(t)\|^{2})^{2}}.\]
By the Schwartz inequality for the second term, we obtain
\[ K(t)\geq \frac{\|\dot{\psi}(t)\|^{2}}{1-\|\psi(t)\|^{2}}
- \frac{  \|\dot{\psi}(t)\|^{2}    \|\psi(t)\|^{2}   }{(1-\|\psi(t)\|^{2})^{2}}
= \frac{\|\dot{\psi}(t)\|^{2}}{(1-\|\psi(t)\|^{2})^{2}}.\]
Hence $K(t)=\frac{\|\dot{\psi}(t)\|^{2}}{(1-\|\psi(t)\|^{2})^{2}}$
if and only if
there exists a curve $\alpha :[0,1]\to {\bf C}$ such that
$\dot{\psi}(t)=\alpha(t)\psi(t)$.
The solution of this equation is given by
$\psi(t)=C{\rm exp}(\int_{0}^{t}\alpha(t)dt)\cdot z$
where $C$ is an integral constant.
Since the length of geodesic is invariant for reparametrization
and the initial condition of $\psi$ 
is $\phi(0)=0$, we find $\psi(t)=t\cdot z$.
Hence the first statement is verified.
Applying $K(t)=\frac{\|z\|^{2}}{(1-t^{2}\|z\|^{2})^{2}}$ to (\ref{eqn:ie}), 
we obtain the second statement. \qedh

\noindent
{\bf Proof of Theorem  \ref{Thm:tr}.}
%
%
We calculate the distance on $B_{{\cal H}}$ according to 
{\it ``Exercises and Further Results,
G. The Hyperbolic Plane"} in $\S$ 1 of \cite{hel}
and $\S$ 2 of \cite{Kobayashi}.
If $u=0$ or $v=0$, then 
it is given by Lemma \ref{lem:16} and in this case, 
the statement is true. 
Assume $u\ne 0$ and $v\ne 0$. 
If $v=k\cdot u$ for some $k\in {\bf R}$,
then $u, v,0$ are on a same geodesic by Lemma \ref{lem:16}.
Hence we assume that $u$ and  $v$ are ${\bf R}$-linearly independent.

Let $e_{1}\equiv u/\|u\|$, $x={\rm Re}\langle e_{1}|v\rangle $, 
$e_{2}\equiv\frac{v-xe_{1}}{\|v-xe_{1}\|}$,
$y\equiv {\rm Re}\langle e_{2}|v\rangle $. Then
 $x=p/\|u\|$, $y=q/\|u\|$, $p^{2}+q^{2}=\|u\|^{2}\|v\|^{2}$ and 
\[v=xe_{1}+ye_{2}=(pe_{1}+qe_{2})/\|u\|.\]
Let
$T=\left(
\begin{array}{cc}
A & m_{u}u \\
 m_{u}\bar{u} & m_{u}
\end{array}
\right )\in U_{1}({\cal H})$
where $m_{u}\equiv (1-\|u\|^{2})^{-1/2}$.	
Then we see that $\phi_{T}(0)=u$
for any $A\in {\cal L}({\cal H})$ satisfying 
(\ref{eqn:conditions}) for $x=y=m_{u}u$ and $a=m_{u}$.

From Lemma \ref{lem:16},
the distance $d(0,w)$ between $0$ and $w\equiv \phi_{T^{-1}}(v)$ is given by
\[d(0,w)=\frac{1}{2}\log \frac{1+\|w\|}{1-\|w\|}. \]
Because $\phi_{T^{-1}}$ is also an isometry, $d(u,v)=d(0,w)$. 
By definition of $w$ and $v$,
\[w=\{ (p/\|u\|-\|u\|)e_{1} + (qm_{u}^{-1}/\|u\|)e_{2}\} (1-p)^{-1} .\]
From this,
we obtain $\|w\|^{2}(1-p)^{2}= \|u-v\|^{2}-q^{2}$.
Hence the statement holds.
\qedh

\noindent
{\bf Proof of Lemma \ref{lem:grad}.}
%
%
(i) 
Let $W\equiv {\rm grad}f$.
By definition, 
$g_{z}(W_{z}, \bar{v})=\bar{\partial}_{z}f(\bar{v})$.
The L.H.S. of this is
\[ k_{z}\langle v|W_{z}\rangle 
+k_{z}^{2}\langle v|z\rangle \langle z|W_{z}\rangle
=\langle  v|(k_{z}W_{z}
+k_{z}^{2}\langle z|W_{z}\rangle z)\rangle \]
from (\ref{eqn:metric}).
From this,
%
%
\begin{equation}\label{eqn:42}
(\bar{\partial}_{z}f)^{*}=k_{z}W_{z}+k_{z}^{2}\langle z|W_{z}\rangle z.
\end{equation}
Hence
%
%
\begin{equation}\label{eqn:43}
 \langle z|(\bar{\partial}_{z}f)^{*}\rangle =(k_{z}+k_{z}^{2}\|z\|^{2})
\langle z|W_{z}\rangle  =k_{z}^{2}\langle z|W_{z}\rangle. 
\end{equation}
By inserting (\ref{eqn:43}) to (\ref{eqn:42}),
$W_{z}=k_{z}^{-1}\{(\bar{\partial}_{z}f)^{*}
-\langle z|(\bar{\partial}_{z}f)^{*}\rangle z\}$.

\noindent
(ii) 
Let $W\equiv {\rm sgrad}f$.
Then
$\omega_{z}(W_{z},\bar{u})=\bar{\partial}_{z}f(\bar{u})$.
The L.H.S. of this is 
\[g_{z}(J_{z}W_{z},\bar{u})
=\sqrt{-1}(k_{z}\langle u|W_{z}\rangle +k_{z}^{2}
\langle u|z\rangle \langle z|W_{z}\rangle ).\]
Hence we obtain that
$(\bar{\partial}_{z}f)^{*}=\sqrt{-1}\{k_{z}W_{z}
+k_{z}^{2}\langle z|W_{z}\rangle z\}$.
As same way about (i), we can verify the statement of (ii).
\qedh

%
%
\begin{lem}
\label{lem:a1}
Let $E_{1}$ be the projection from $\widetilde{{\cal H}}$ to ${\cal H}$ and 
$e_{2}\equiv (0,1)\in \widetilde{{\cal H}}$.
Then the following holds:
%
%
\begin{eqnarray}
(\partial_{z}f_{C})(u)& =&
k_{z}\langle \hat{z}|Cu\rangle+k_{z}^{2}
\langle\hat{z}|C\hat{z}\rangle\langle z|u\rangle,
\label{eqn:ew1}\\
(\partial_{z}f_{C})^{*}& =&
k_{z}E_{1}C^{*}\hat{z}+k_{z}^{2}\langle\hat{z}|C^{*}\hat{z}\rangle z, 
\label{eqn:ew2}\\
{\rm grad}_{z}f_{C}&=&E_{1}C\hat{z}+\langle e_{2}|C\hat{z}\rangle z 
\label{eqn:ew3}
\end{eqnarray}
where $(\partial_{z}f_{C})^{*}$ is 
the holomorphic tangent vector at $z$ defined by
$\langle (\partial_{z}f_{C})^{*}|X\rangle=(\partial_{z}f_{C})(X)$
for $X\in T_{z}B_{{\cal H}}$.
\end{lem}
%
%
\pr
We show only the later two parts.
From (\ref{eqn:ew1}),
\[
\begin{array}{ll}
 (\partial_{z}f_{C})(u)&=
\langle k_{z}C^{*}\hat{z}|u\rangle +\langle \,
k_{z}^{2}\overline{\langle \hat{z}|C\hat{z}\rangle}
z\,|u\rangle 
 \\
&=\langle \{k_{z}C^{*}\hat{z}+k_{z}^{2}
\langle \hat{z}|C^{*}\hat{z}\rangle z\}|u\rangle \\
&=\langle \{k_{z}C^{*}\hat{z}+k_{z}^{2}\langle \hat{z}|C^{*}\hat{z}\rangle z\}|E_{1}u\rangle \\
&=\langle \{k_{z}E_{1}C^{*}\hat{z}+k_{z}^{2}\langle \hat{z}|C^{*}\hat{z}\rangle z\}|u\rangle .
\end{array}
\]
By comparing both side of this equation, we obtain (\ref{eqn:ew2}). 
By inserting (\ref{eqn:ew2}) to Lemma \ref{lem:grad}, we find (\ref{eqn:ew3}).
\qedh

%

\end{document}